# Nearly-integrable perturbations of the Lagrange top: applications of KAM-theory[*]

H. W. Broer[1], H. Hanßmann[2,3], J. Hoo[1] and V. Naudot[1]

*Groningen University, RWTH Aachen and Groningen University*

**Abstract:** Motivated by the Lagrange top coupled to an oscillator, we consider the quasi-periodic Hamiltonian Hopf bifurcation. To this end, we develop the normal linear stability theory of an invariant torus with a generic (i.e., non-semisimple) normal $1 : -1$ resonance. This theory guarantees the persistence of the invariant torus in the Diophantine case and makes possible a further quasi-periodic normal form, necessary for investigation of the non-linear dynamics. As a consequence, we find Cantor families of invariant isotropic tori of all dimensions suggested by the integrable approximation.

## 1. The Lagrange top

The Lagrange top is an axially symmetric rigid body in a three dimensional space, subject to a constant gravitational field such that the base point of the body-symmetry (or figure) axis is fixed in space, see Figure 1. Mathematically speaking, this is a Hamiltonian system on the tangent bundle $TSO(3)$ of the rotation group $SO(3)$ with the symplectic 2-form $\sigma$. This $\sigma$ is the pull-back of the canonical 2-form on the co-tangent bundle $T^*SO(3)$ by the bundle isomorphism $\tilde{\kappa} : TSO(3) \to T^*SO(3)$ induced by a non-degenerate left-invariant metric $\kappa$ on $SO(3)$, where $\tilde{\kappa}(v) = \kappa(v, \cdot)$. The Hamiltonian function $H$ of the Lagrange top is obtained as the sum of potential and kinetic energy. In the following, we identify the tangent bundle $TSO(3)$ with the product $M = SO(3) \times \mathfrak{so}(3)$ via the map $v_Q \in TSO(3) \mapsto (Q, T_{\mathrm{Id}}L_Q^{-1}v) \in SO(3) \times \mathfrak{so}(3)$, where $\mathfrak{so}(3) = T_{\mathrm{Id}}SO(3)$ and $L_Q$ denotes left-translation by $Q \in SO(3)$. We assume that the gravitational force points vertically downwards. Then, the Lagrange top has two rotational symmetries: rotations about the figure axis and the vertical axis $e_3$. We let $\mathcal{S} \subset SO(3)$ denote the subgroup of rotations preserving the vertical axis $e_3$. Then, for a suitable choice of the space coordinate system $(e_1, e_2, e_3)$, the two symmetries correspond to a symplectic right action $\Phi^r$ and to a symplectic left action $\Phi^l$ of the Lie subgroup $\mathcal{S}$ on $M$. By the Noether Theorem [1, 4], these Hamiltonian symmetries give rise to

---

[*]Work supported by grant MB-G-b of the Dutch FOM program Mathematical Physics and by grant HPRN-CT-2000-00113 of the European Community funding for the Research and Training Network MASIE.

[1]Department of Mathematics, University of Groningen, PO Box 800, 9700 AV Groningen, The Netherlands, e-mail: `h.w.broer@math.rug.nl` e-mail: `j.hoo@math.rug.nl` e-mail: `v.naudot@math.rug.nl`

[2]Institut für Reine & Angewandte Mathematik, RWTH Aachen, 52056 Aachen, Germany.

[3]Present address: Mathematisch Instituut, Universiteit Utrecht, Postbus 80010, 3508 TA Utrecht, The Netherlands, e-mail: `heinz.hanssmann@math.uu.nl`

*AMS 2000 subject classifications:* primary 37J40; secondary 70H08.

*Keywords and phrases:* KAM theory, quasi-periodic Hamiltonian Hopf bifurcation, singular foliation, the Lagrange top, gyroscopic stabilization.

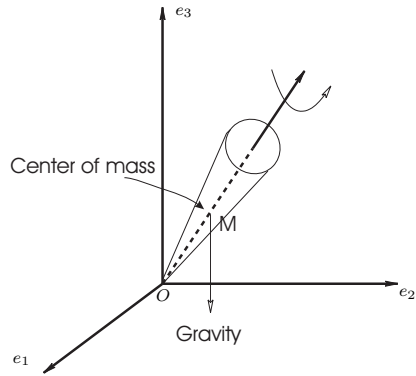

Fig 1. *The Lagrange top in the vertical gravitational force field.*

integrals $\mathcal{M}^r$ and $\mathcal{M}^l$ of the Hamiltonian $H$: the angular momenta along the figure axis and along the vertical axis, respectively. These integrals induce the so-called *energy-momentum map* $\mathcal{EM} := (\mathcal{M}^r, \mathcal{M}^l, H) : M \to \mathbb{R}^3$.

The inverse images of the map $\mathcal{EM}$ divide the phase space $M$ into invariant sets of the Hamiltonian system $X_H$ associated with $H$: for a regular value $m = (a, b, h) \in \mathbb{R}^3$, the set $\mathcal{EM}^{-1}(m) \subset M$ is an $X_H$-invariant 3-torus; for a critical value $m$, this set is a 'pinched' 3-torus, a 2-torus or a circle. This division is a singular foliation of the phase space by $X_H$-invariant tori. This foliation gives rise to a stratification of the parameter space: the $(a, b, h)$-space is split into different regions according to the dimension of the tori $\mathcal{EM}^{-1}(a, b, h)$. To describe the stratification, one applies a regular reduction [1, 28] by the right symmetry $\Phi^r$ to the three-degrees-of-freedom Hamiltonian $H$ as follows. For a fixed value $a$, we deduce from $H$ a two-degrees-of-freedom Hamiltonian $H_a$ on the four-dimensional orbit space $M_a = (\mathcal{M}^r)^{-1}(a)/\mathcal{S}$ under the $\Phi^r$-action. The reduced phase space $M_a$ can be identified with the four-dimensional submanifold $\mathcal{R}_a = \{(u, v) \in \mathbb{R}^3 \times \mathbb{R}^3 : u \cdot u = 1, u \cdot v = a\}$, compare with [13, 16]. This manifold $\mathcal{R}_a$ inherits the symplectic 2-form $\omega_a$ given by

$$\omega_a(u, v)\big((x, y), (p, q)\big) = \hat{x} \cdot \hat{q} - \hat{y} \cdot \hat{p} + v \cdot (\hat{v} \times \hat{p}), \tag{1.1}$$

where $(x, y) = (\hat{x} \times u, \hat{x} \times v + \hat{y}) \in T_{(u,v)}\mathcal{R}_a$, $(p, q) = (\hat{p} \times u, \hat{p} \times v + \hat{q}) \in T_{(u,v)}\mathcal{R}_a$. Here $\cdot$ and $\times$ denote the standard inner and cross product of $\mathbb{R}^3$, respectively. The reduced Hamiltonian $H_a : \mathcal{R}_a \to \mathbb{R}$ obtains the form

$$H_a(u, v) = \frac{1}{2} v \cdot v + c u_3 + \rho a^2,$$

where $c > 0$ and $\rho \in \mathbb{R}$. Observe that $P_a = (0, 0, 1, 0, 0, a) \in \mathcal{R}_a$ is an isolated equilibrium of $H_a$. For each $a$, the point $P_a$ — a relative equilibrium of the full system $H$ — corresponds to a periodic solution of $H$, namely a rotation about the vertical axis.

It is known [13, 17] that a stability transition of the rotations $P_a$ occurs as the angular momentum value $a$ passes through the critical value $a = a_0 = 2\sqrt{c}$. Physically, this transition is referred to as *gyroscopic stabilization* of the Lagrange top. A mathematical explanation for the stabilization is that Floquet matrices $\Omega_a$ of the periodic orbits $P_a$ changes from hyperbolic into elliptic as $a$ passes through $a_0$, see Figure 2. In fact, the Lagrange top undergoes a (non-linear) *Hamiltonian Hopf bifurcation* [13, 17]. A brief discussion of this bifurcation is given in Section 2.1 below; for an extensive treatment see [29]. Following [13, 17, 29, 30], the local stratification of the parameter space — associated with the singular foliation by invariant tori — near the bifurcation point is described by a piece of the swallowtail catastrophe set from singularity theory [37] as follows. The one-dimensional singular part of this surface is the stratum associated with the periodic solutions, the regular part forms the stratum of 2-tori and the open region above the surface is the stratum of Lagrangian 3-tori, compare with Figure 3.

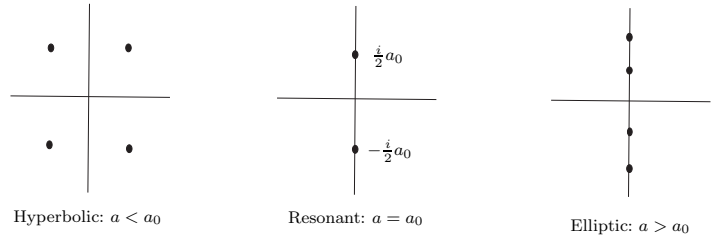

FIG 2. *Eigenvalue configuration of Floquet matrix $\Omega_a$ as $a$ passes through $a_0$.*

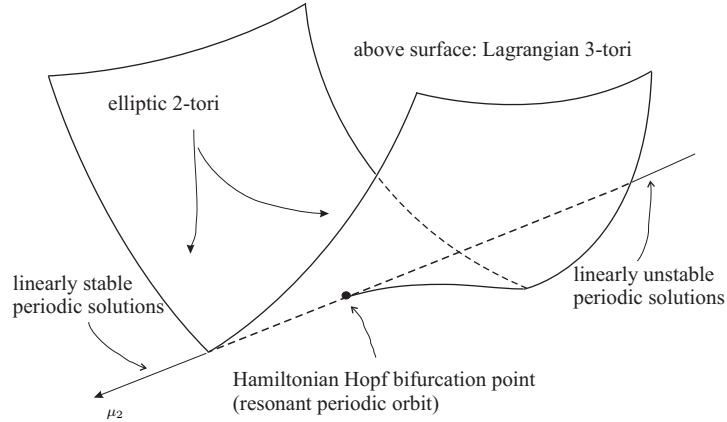

FIG 3. *Sketch of the local stratification by invariant tori near the Hamiltonian Hopf bifurcation point in the $(a, b, h)$-space. After a suitable reparametrization, the surface is a piece of the swallowtail catastrophe set [17, 29, 37]. The parameter $\mu_2$ is given by $\frac{a^2}{4} - c$.*

We are interested in a perturbation problem where the Lagrange top is weakly coupled to an oscillator with multiple frequencies, e.g., the base point of the top is coupled to a vertically vibrating table-surface by a massless spring, see Figure 4. In this example, the spring constant, say $\varepsilon$, plays the role of the perturbation parameter.

More generally, we consider the Hamiltonian perturbation $H_\varepsilon$ of the form

$$H_\varepsilon = H + G + \varepsilon F, \quad (1.2)$$

where $H$ and $G$ are the Hamiltonians of the Lagrange top and of the quasi-periodic oscillator, respectively. Here the function $F$ depends on the coupling between the top and the oscillator. We assume that the quasi-periodic oscillator has $n \geq 1$ frequencies and is Liouville-integrable [1, 4]. Then, the unperturbed integrable Hamiltonian $H_0 = H + G$ contains invariant $(n+1)$-, $(n+2)$- and $(n+3)$-tori. Near gyroscopic stabilization of the Lagrange top, this Hamiltonian $H_0$ gives a similar (local) stratification by tori as sketched in Figure 3 in a parameter space, but with invariant $(n+1)$-, $(n+2)$- and $(n+3)$-tori, for details see [7]. Our concern is with the fate of these invariant tori for small but non-zero $\varepsilon$.

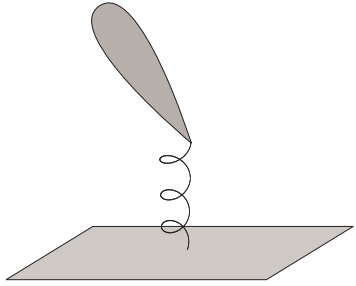

FIG 4. *The Lagrange top coupled with a vibrating table-surface by a spring.*

By KAM theory [2, 10, 12, 27, 33, 35, 39], the 'majority' of the invariant tori from the local stratification survives small perturbations. They form Whitney-smooth Cantor families, parametrized over domains with positive measure. The purpose of the present paper is to describe the persisting families of tori in terms of a quasi-periodic Hamiltonian Hopf bifurcation [7, 10]. From [13, 17] it is known that the Lagrangian torus bundle in the unperturbed Lagrangian top contains non-trivial monodromy. More precisely, we consider the local stratification by tori as sketched in Figure 3. Let $D$ be a punctured disk in the stratum of Lagrangian 3-tori which transversally intersects the 'thread' (the 1-dimensional curve associated with unstable periodic solutions). Then, the Lagrangian torus bundle with the boundary $\partial D \cong \mathbb{S}^1$ as the base space is non-trivial. From this we conclude that the Lagrangian $(n+2)$-torus bundle of the unperturbed Hamiltonian $H_0$ has non-trivial monodromy. We like to mention that, in view of the global KAM theory [6], there exists a proper extension of the non-trivial monodromy in the integrable Lagrangian torus bundle, to the nearly-integrable one. In this sense, we may say that the non-trivial monodromy, that goes with the Hamiltonian Hopf bifurcation and is centered at the thread, survives a small non-integrable perturbation.

## 2. Hamiltonian Hopf bifurcations of equilibria, periodic and quasi-periodic solutions

The above persistence problem is part of a more general study of a quasi-periodic Hamiltonian Hopf bifurcation to be discussed in Section 2.2 below. This quasi-periodic bifurcation can be considered as a natural extention of the Hamiltonian Hopf bifurcation of equilibria [29]. Let us first recall certain facts about the latter.

### 2.1. The Hamiltonian Hopf bifurcation of equilibria

We consider a two-degrees-of-freedom Hamiltonian on $\mathbb{R}^4 = \{z_1, z_2, z_3, z_4\}$ with the standard symplectic 2-form $dz_1 \wedge dz_3 + dz_2 \wedge dz_4$. Let $\hat{H} = \hat{H}_\mu$ be a family of Hamiltonian functions with the origin as an equilibrium. Assume that for $\mu = \mu_0$, the linearized Hamiltonian system (at the origin) has a double pair of purely imaginary eigenvalues with a non-trivial nilpotent part.[1] Moreover, as $\mu$ passes through the value $\mu_0$, the eigenvalues of the linear part (at the origin) behave as follows: as $\mu$ increases a complex quartet moves towards the imaginary axis, meeting there for $\mu = \mu_0$ and splitting into two distinct purely imaginary pairs for $\mu > \mu_0$, compare with Figure 2. This bifurcation is referred to as a *generic Hamiltonian Hopf bifurcation*, provided that a certain generic condition on the higher order terms is met.

---

[1] In this case, we say that the Hamiltonian $\hat{H}$ is in generic or (non-semisimple) $1:-1$ resonance at the origin for $\mu = \mu_0$. Notice that the linear part at the origin is indefinite.

More precisely, we normalize the Hamiltonian $\hat{H}$ with respect to the quadratic part of $\hat{H}$. This gives us the Birkhoff normal form

$$\begin{aligned}\hat{H} &= (\nu_1(\mu) + \lambda_0)S + N + \nu_2(\mu)M + b_1(\mu)M^2 \\ &\quad + b_2(\mu)SM + b_3(\mu)S^2 \ + \ \text{h.o.t.},\end{aligned} \qquad (2.1)$$

where $M = \frac{1}{2}(z_1^2 + z_2^2)$, $S = z_1 z_4 - z_2 z_3$ and $N = \frac{1}{2}(z_3^2 + z_4^2)$, see [24, 29, 30]. Here $\lambda_0 \neq 0$, and $\nu_1(0) = 0 = \nu_2(0)$. The generic condition now requires that $\frac{\partial \nu_2}{\partial \mu}(\mu_0) \neq 0$ and $b_1(\mu_0) \neq 0$. We focus on the *supercritical* case where $b_1(\mu_0) > 0$. We may assume that $\hat{H}$ is invariant under the $\mathbb{S}^1$-action generated by the flow of the semi-simple quadratic part $S$, since this symmetry can be pushed through the normal form (2.1) up to an arbitrary order [29]. Then, the periodic solutions of $\hat{H}$ are given by the singularities of the energy-momentum map $(\hat{H}, S)$. Following [29], this map has a (singularity theoretical) normal form $(G, S)$ where $G = N + \delta(\mu)M + M^2$. Indeed, the maps $(\hat{H}, S)$ and $(G, S)$ are locally left-right equivalent by an $\mathbb{S}^1$-equivariant origin-preserving diffeomorphism on $\mathbb{R}^4$ and an origin-preserving diffeomorphism on $\mathbb{R}^2$. As a result, the set of critical values $\mathcal{C}$ of $(G, S)$, considered as the graph over the parameter $\mu$, is diffeomorphic to that of $(\hat{H}, S)$, which is a piece of swallowtail surface [37] in the $(\delta, S, G)$-space. This critical set $\mathcal{C}$ determines the local stratification by leaves of $(G, S)$ near the bifurcation point in the parameter space: strata of equilibria, periodic orbits and 2-tori. This stratification corresponds to the situation as sketched in Figure 3, when replacing periodic solutions, 2-tori and 3-tori by equilibria, periodic orbits and 2-tori, respectively.

**Remarks 2.1.**

1. The above discussion also applies for the situation where $\mu$ is a multi-parameter [30].
2. In the case where the Hamiltonian $\hat{H}$ is not $\mathbb{S}^1$-symmetric, one first applies a Liapounov-Schmidt reduction to obtain an $\mathbb{S}^1$-symmetric Hamiltonian $H$ which has the same normal form as $\hat{H}$ up to arbitrary order [29, 30, 42]. This reduction relates the periodic solutions of $\hat{H}$ to that of $H$ in a diffeomorphic way.

### 2.2. The Hamiltonian Hopf bifurcation of (quasi-) periodic solutions

Motivated by the persistence problem of the Lagrange top coupled to an oscillator, see Section 1, we consider the quasi-periodic dynamics of Hamiltonian systems with more degrees of freedom near an invariant resonant torus, where our main interest is with the normal $1 : -1$ resonance. More precisely, our phase space $M$ is given by $\mathbb{T}^m \times \mathbb{R}^m \times \mathbb{R}^4 = \{x, y, z\}$ with symplectic 2-form $\sigma = \sum dx_i \wedge dy_i + \sum dz_j \wedge dz_{j+2}$, where $\mathbb{T}^m = \mathbb{R}^m/(2\pi\mathbb{Z})^m$. This space $M$ admits a free $\mathbb{T}^m$-action given by $(\theta, (x, y, z)) \in \mathbb{T}^m \times M \mapsto (\theta + x, y, z) \in M$. A Hamiltonian function is said to be $\mathbb{T}^m$-*symmetric* or *integrable* if it is invariant under this $\mathbb{T}^m$-action, compare with

[12, 26]. We consider a $p$-parameter family of integrable Hamiltonian functions

$$H(x,y,z;\nu) = \langle \omega(\nu), y \rangle + \frac{1}{2}\langle Jz, \Omega(\nu)z \rangle + \text{h.o.t.}, \tag{2.2}$$

where $\nu \in \mathbb{R}^p$ is the parameter, $\omega(\nu) \in \mathbb{R}^n$, $\Omega(\nu) \in \mathfrak{sp}(4, \mathbb{R})$ and $J$ is the standard symplectic $4 \times 4$-matrix. By $\mathbb{T}^m$-symmetry all higher order terms are $x$-independent. Then the union $T = \bigcup_\nu T_\nu \subseteq M \times \mathbb{R}^p$, where $T_\nu = \{(x, y, z, \nu) : (y, z) = (0, 0)\}$, is a $p$-parameter family of invariant $m$-tori of $H$ parametrized by $\nu$. Let the torus $T_{\nu_0}$ be in *generic normal* $1 : -1$ *resonance*, meaning that the Floquet matrix $\Omega(\nu_0)$ at $T_{\nu_0}$ has a double pair of purely imaginary eigenvalues with a non-trivial nilpotent part. Note that for $m = 0$ we arrive at the setting of the equilibria case [29], compare with Section 2.1. For $m = 1$, the invariant submanifolds $T_\nu$ are periodic solutions of the Hamiltonian $H$. For this reason, we speak of the periodic case. Similarly, we refer to our present setting with $m \geq 2$ as the quasi-periodic case.

Let us consider the local (quasi-)periodic dynamics of the family $H$ with $m \geq 1$ near the resonant torus $T_{\nu_0}$. To this end, we first reduce the Hamiltonian $H$ by the $\mathbb{T}^m$-symmetry as follows. For a small fixed value $\alpha \in \mathbb{R}^m$, we obtain from the Hamiltonian $H$ a two-degrees-of-freedom $H_\alpha$ defined on the space $M_\alpha = (M/\mathbb{T}^m) \cap \{y = \alpha\}$, where $M/\mathbb{T}^m$ denotes the orbit space of the $\mathbb{T}^m$-action. More explicitly, we obtain $H_\alpha(z;\nu) = H(x, \alpha, z; \nu)$, by identifying $M_\alpha$ with $\mathbb{R}^4 = \{z_1, z_2, z_3, z_4\}$.

We may assume that the reduced Hamiltonian $H_\alpha$ is in normal form with respect to the $\mathbb{S}^1$-symmetry generated by the semi-simple part of the quadratic term $\frac{1}{2}\langle Jz, \Omega(\nu_0)z \rangle$ [29]. Then, the full Hamiltonian $H$ is invariant under this circle-action. We say that the full Hamiltonian $H$ undergoes a generic *periodic* ($m = 1$) or *quasi-periodic* ($m \geq 2$) *Hamiltonian Hopf bifurcation* at $\nu = \nu_0$, if the reduced two-degrees-of-freedom system $H_\alpha$ has a generic Hamiltonian Hopf bifurcation, see Section 2 and Remark 2.1. As before we restrict to the supercritical case. In the following, we focus on the quasi-periodic case where $m \geq 2$. For a treatment of the periodic case see [34]. Applying the local analysis of Section 2.1 to the reduced Hamiltonian $H_\alpha$, we conclude that the local torus foliation of $M$ near the resonant torus $T_{\nu_0}$ defines a local stratification in a parameter space by $m$-, $(m+1)$- and $(m+2)$-tori of the full system $H$. This stratification is sketched in Figure 5, compare with Figure 3. Notice that the resonant torus $T_{\nu_0}$ corresponds to the quasi-periodic Hamiltonian Hopf bifurcation point.

Our main goal is to investigate the persistence of the invariant $m$-, $(m+1)$- and $(m+2)$-tori from the local stratification as in Figure 5(a), when the integrable Hamiltonian $H$ is perturbed into a nearly-integrable (i.e., not necessarily $\mathbb{T}^m$-symmetric) one. Observe that in the example of the Lagrange top coupled with a quasi-periodic oscillator one has $m = n + 1$, see Section 1. In the sequel we examine the persistence of the $p$-parameter family $T = \bigcup_\nu T_\nu$ of invariant $m$-tori. This family consists of elliptic, hyperbolic tori and the resonant torus $T_{\nu_0}$, compare with Figure 5(a). The 'standard' KAM theory [2, 12, 26, 31, 32, 33, 35, 36, 39] yields persistence only for subfamilies of the family $T$, containing elliptic or hyperbolic tori. The problem is that the resonant torus $T_{\nu_0}$ gives rise to multiple Floquet exponents, compare with Figure 2. To deal with this problem, we develop a normal linear stability theorem [10], as an extension of the 'standard' KAM theory. The persistence of the $(m+1)$- and $(m+2)$-tori will be discussed in Section 4.

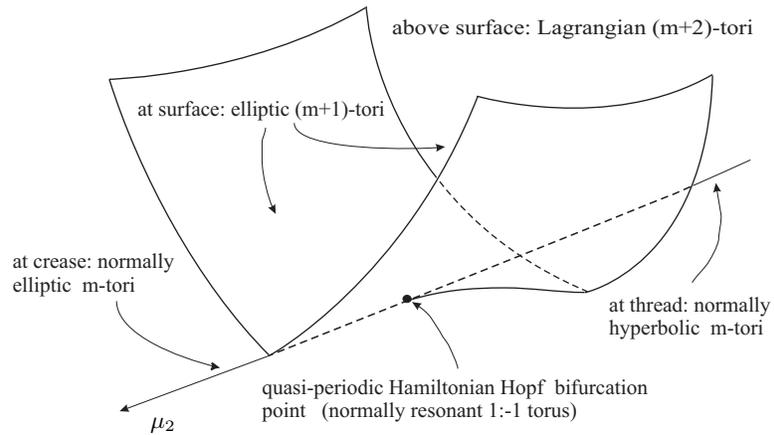

(a) Unperturbed situation

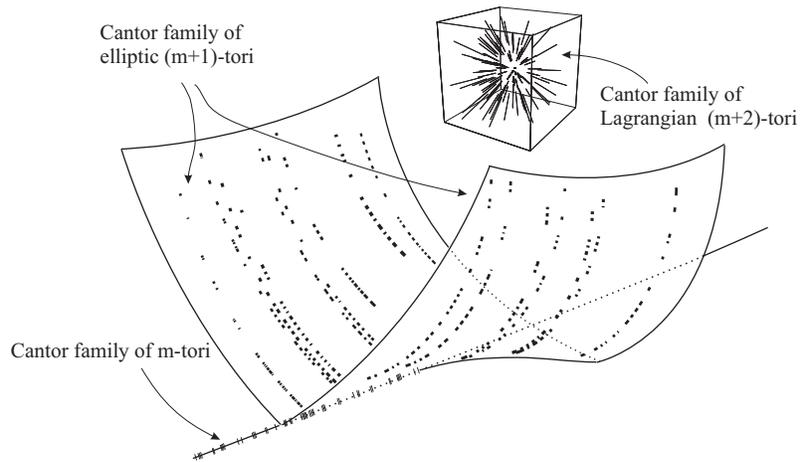

(b) Perturbed situation

FIG 5. (a) *Singular foliation by invariant tori of unperturbed integrable Hamiltonian $H$ near the resonant torus in the supercritical case;* (b) *sketch of the Cantor families of surviving Diophantine invariant tori in the singular foliation of the perturbed, nearly-integrable Hamiltonian $\widetilde{H}$.*

## 3. The quasi-periodic Hamiltonian Hopf bifurcation: persistence of the $m$-tori

In this section we develop the normal linear stability theory, needed for the persistent $m$-tori [10]. This theory also allows us to obtain a further quasi-periodic normal form of perturbations, necessary for investigation of the non-linear nearly-integrable dynamics.

### 3.1. Normal linear stability: a part of KAM theory

Let us consider a general $p$-parameter family $K = K_\nu$ of integrable real-analytic Hamiltonian functions on the phase space $M = \mathbb{T}^m \times \mathbb{R}^m \times \mathbb{R}^{2q} = \{x, y, z\}$ with the symplectic 2-form $\sum_{i=1}^{m} dx_i \wedge dy_i + \sum_{j=1}^{q} dz_j \wedge dz_{j+2}$. We assume that the Hamiltonian vector field $X = X_K$ associated with $K$ is of the form

$$X_\nu(x, y, z) = [\omega(\nu) + O(|y|, |z|)]\frac{\partial}{\partial x} + [\Omega(\nu)z + O(|y|, |z|^2)]\frac{\partial}{\partial z}, \qquad (3.1)$$

where $\omega(\nu) \in \mathbb{R}^m$ and $\Omega(\nu) \in \mathfrak{sp}(2q, \mathbb{R})$. Then, the torus family $T = \cup_\nu T_\nu$, where

$$T_\nu = \{((x, y, z), \nu) \in M \times \mathbb{R}^p : (y, z) = (0, 0)\},$$

is an $X_K$-invariant submanifold. Note that the integrable Hamiltonian $H$ given by (2.2) is a special form of the family $K$. A typical KAM-stability question is concerned with the persistence of the invariant submanifold $T$ under small perturbation of the integrable Hamiltonian family $K$. We refer to $\Omega(\nu)$ as the Floquet (or normal) matrix of the torus $T_\nu$. The 'standard' KAM theory [2, 12, 26, 27, 31, 32, 33, 35, 36, 39] asserts that the 'majority' of these invariant tori survives small perturbations, provided that the unperturbed family satisfies the following conditions:

(a) the Floquet exponents of the torus $T_{\nu_0}$ (i.e., the eigenvalues of $\Omega(\nu_0)$) are simple;
(b) the matrix $\Omega(\nu_0)$ is non-singular;
(c) the product map $\omega \times \Omega : \mathbb{R}^p \to \mathbb{R}^m \times \mathfrak{sp}(2q, \mathbb{R})$ is transversal to the submanifold $\{\omega(\nu_0)\} \times \text{Orbit}\, \Omega(\nu_0)$, where $\text{Orbit}\, \Omega(\nu_0)$ denotes the similarity class of $\Omega(\nu_0)$ by the linear symplectic group;
(d) the internal and normal frequencies satisfy Diophantine conditions.

Let us be more specific about these assumptions. The internal frequencies of the invariant torus $T_\nu$ are given by $\omega(\nu) = (\omega_1(\nu), \ldots, \omega_m(\nu))$, and the normal frequencies consist of the positive imaginary parts of the eigenvalues of $\Omega(\nu)$. Conditions (a) and (b) require that the eigenvalues of $\Omega_0 = \Omega(\nu_0)$ have distinct non-zero eigenvalues. Condition (c) means that the map $\omega$ is submersive at $\nu = \nu_0$ and the map $\Omega$ is a versal unfolding of $\Omega_0$ (in the sense of [3, 22]) simultaneously. The Hamiltonian $K$ is said to be non-degenerate at the torus $T_{\nu_0}$, if it meets conditions (b) and (c).

**Remark 3.1.** The non-degeneracy has the following geometrical interpretation: the normal linear part $NX_\nu = \omega(\nu)\frac{\partial}{\partial x} + \Omega(\nu)z\frac{\partial}{\partial z}$ is *transversal* to the conjugacy class of $NX_0$ within the space of normally affine Hamiltonian vector fields as in [3, 10, 22].

Denoted by $\omega^N(\nu) = (\omega_1^N(\nu), \ldots, \omega_r^N(\nu))$ — called the normal frequencies of the torus $T_\nu$ — condition (d) is formulated as follows: for a constant $\tau > m - 1$ and a parameter $\gamma > 0$, we have that

$$\left|\langle \omega(\nu), k \rangle + \langle \omega^N(\nu), \ell \rangle\right| \geq \gamma |k|^{-\tau}, \qquad (3.2)$$

for all $k \in \mathbb{Z}^m \setminus \{0\}$ and for all $\ell \in \mathbb{Z}^2$, where $|\ell_1| + \cdots + |\ell_r| \leq 2$. Due to the simpleness assumption (a), the number $r$ of the normal frequencies is independent of $\nu$, for $\nu$ sufficiently close to $\nu_0$, compare with [12, 26]. The map $\mathcal{F} : \nu \mapsto (\omega(\nu), \omega^N(\nu))$ is called the frequency map. We denote by $\Gamma_{\tau,\gamma}(U)$ the set of parameters $\nu \in U$ such that $\mathcal{F}(\nu)$ satisfy the non-resonant condition (3.2). Also we need the subset $\Gamma_{\tau,\gamma}(U')$, where the set $U' \subset U$ is given by

$$U' = \left\{ \nu \in U : \operatorname{dist.}\big((\omega(\nu), \omega^N(\nu)), \partial \mathcal{F}(U)\big) > \gamma \right\}. \tag{3.3}$$

Observe that for $\gamma$ sufficiently small, the set $U'$ is still an non-empty open neighbourhood of $\nu_0$ and that $\Gamma_{\tau,\gamma}(U') \subset U' \subset U$ contains a 'Cantor set' of Diophantine frequencies with positive measure.

**Remarks 3.2.**

- Condition (a) ensures that the Floquet matrices $\Omega(\nu)$ are semi-simple and that the normal frequencies (after a suitable reparametrization) depend on parameters in an affine way.
- Let $\Lambda$ be the set of $(\omega, \omega^N) \in \mathbb{R}^m \times \mathbb{R}^r$ satisfying the Diophantine conditions (3.2). Then $\Lambda$ is a nowhere dense, uncountable union of closed half lines. The intersection $\Lambda \cap \mathbb{S}^{m+r-1}$ with the unit sphere of $\mathbb{R}^m \times \mathbb{R}^r$ is a closed set, which by Cantor-Bendixson theorem [25] is the union of a perfect set $\mathcal{P}$ and a countable set. Note that the complement of $\Lambda \cap \mathbb{S}^{m+r-1}$ contains the dense set of resonant vectors $(\omega, \omega^N)$. Since all points of $\Lambda$ are separated by the resonant hyperplanes, this perfect set $\mathcal{P}$ is totally disconnected and hence a Cantor set. In $\mathbb{S}^{m+r-1}$ this Cantor set tends to full Lebesgue measure as $\gamma \downarrow 0$. We refer to the set of frequencies $(\omega(\nu), \omega^N(\nu))$ satisfying (3.2) as a 'Cantor set' — a foliation of manifolds over a Cantor set.

Though condition (a) is generic,[2] it excludes certain interesting examples like the quasi-periodic Hamiltonian Hopf bifurcation, as it occurs in the example of the Lagrange top coupled to a quasi-periodic oscillator. Indeed, such systems have a Floquet matrix with multiple eigenvalues, compare with Figure 2. To deal with this problem, we have to drop the simpleness assumption. Instead we impose a Hölder condition on the spectra $\operatorname{Spec} \Omega(\nu)$ of the matrices $\Omega(\nu)$ as follows. Let $U \subset \mathbb{R}^p$ be a small neighbourhood of $\nu_0$. Suppose that $\Omega$ has a holomorphic extension to the complex domain

$$U + r_0 = \{\tilde{\nu} \in \mathbb{C}^p : \exists \nu \in U \text{ such that } |\nu - \tilde{\nu}| \leq r_0\} \subset \mathbb{C}^p \tag{3.4}$$

for a certain constant $r_0 > 0$. We say that $\Omega(\nu)$ is $(\theta, r_0)$-Hölder, if there exist positive constants $\theta$ and $L$ such that the following holds: for any $\tilde{\nu} \in U + r_0, \nu \in U$ and for any $\tilde{\lambda} \in \operatorname{Spec} \Omega(\tilde{\nu})$, there exist a $\lambda \in \operatorname{Spec} \Omega(\nu)$ such that

$$\left| \operatorname{Im} \lambda - \operatorname{Im} \tilde{\lambda} \right| \leq L \left| \nu - \tilde{\nu} \right|^\theta, \tag{3.5}$$

where Im denotes the imaginary part. The condition (3.5) holds in particular for matrices with simple eigenvalues. As an extension of the 'standard' KAM theory, we have the following.

**Theorem 3.3 (Normal linear stability [10]).** *Let $K = K_\nu$ be a p-parameter real-analytic family of integrable Hamiltonians with the corresponding Hamiltonian*

---

[2]The set of simple matrices is dense and open in the matrix space $\mathfrak{gl}(2q, \mathbb{R})$.

vector fields given by (3.1). Suppose that $K$ satisfies the non-degeneracy conditions (b) and (c). Also assume that the matrix family $\Omega(\nu)$ is $(\theta, r_0)$-Hölder, see (3.5). Then, for $\gamma > 0$ sufficiently small and for any real-analytic Hamiltonian family $\widetilde{K}$ sufficiently close to $K$ in the compact-open topology on complex analytic extensions, there exists a domain $\mathcal{U}$ around $\nu_0 \in \mathbb{R}^p$ and a map

$$\Phi : M \times \mathcal{U} \to M \times \mathbb{R}^p,$$

defined near the torus $T_{\nu_0}$, such that,

i. $\Phi$ is a $C^\infty$-near-the-identity diffeomorphism onto its image;
ii. The image of the Diophantine tori $V = \bigcup_{\nu \in \Gamma_{\tau,\gamma}(\mathcal{U}')} ( T_\nu \times \{\nu\} )$ under $\Phi$ is $\widetilde{K}$-invariant, and the restriction of $\Phi$ on $V$ conjugates the quasi-periodic motions of $K$ to those of $\widetilde{K}$;
iii. The restriction $\Phi|_V$ is symplectic and preserves the (symplectic) normal linear part[3] $NX = \omega(\nu)\frac{\partial}{\partial x} + \Omega(\nu)z\frac{\partial}{\partial z}$ of the Hamiltonian vector field $X = X_K$ associated with $K$.

We refer to the Diophantine tori $V$ (also its diffeomorphic image $\Phi(V)$) as a *Cantor family of invariant $m$-tori*, as it is parametrized over a 'Cantor set'. The stability Theorem 3.3 includes the cases where the Floquet matrix $\Omega(\nu_0)$ is in (non-semisimple) $1:-1$ resonant, see [10] for details. For the definition of $1:-1$ resonance see Section 2.2. In particular, it is applicable to our persistence problem as formulated in 2.2, regarding the invariant $m$-tori with a normally $1:-1$ resonant torus, compare with Figure 5(a).

**Remark 3.4.** Theorem 3.3, as is generally the case in the 'standard' KAM-theory, requires the invertibility of the Floquet matrix of the central torus $T_{\nu_0}$. We expect that this assumption can be relaxed by using an appropriate transversality condition, compare with Remark 3.1. For recent work in this direction see [9, 43, 44].

### 3.2. Persistence of m-tori

We return to the setting of Section 2.2. Briefly summarizing, we consider a $p$-parameter real-analytic family $H = H_\nu$ of $\mathbb{T}^m$-symmetric (or integrable) Hamiltonian functions on the space $M = \mathbb{T}^m \times \mathbb{R}^m \times \mathbb{R}^4$ given by (2.2), parametrized over a small neighbourhood $U$ of $\nu_0$. This family has an invariant torus family $T = \bigcup_\nu T_\nu$, where

$$T_\nu = \{(x, y, z, \nu) \in M \times \mathbb{R}^p : (y, z) = (0, 0)\}.$$

Moreover, the Hamiltonian $H$ undergoes a supercritical quasi-periodic Hamiltonian Hopf bifurcation at $\nu = \nu_0$. By $\mathbb{T}^m$-symmetry and an application of [29], near the normally resonant torus $T_{\nu_0}$ there is a local singular foliation by invariant hyperbolic and elliptic $m$-tori, elliptic $(m+1)$-tori and Lagrangian $(m+2)$-tori of the Hamiltonian $H$. The local stratification associated with this foliation in a suitable parameter space is a piece of swallowtail and is sketched by Figure 5(a). Presently, we are concerned with the persistence of these $m$-tori from the local foliation under perturbation.

By assumption, the central $m$-torus $T_{\nu_0}$ is (normally) generically $1:-1$ resonant. We also assume that the Hamiltonian $H$ is non-degenerate at the invariant tori $T_{\nu_0}$, see Section 3.1. As a direct consequence of this non-degeneracy, we need (at least) $m + 2$ parameters, that is, $p \geq m + 2$. Let $X = X_\nu$ be the family of Hamiltonian vector fields corresponding to the Hamiltonian family $H = H_\nu$. By the Inverse Function Theorem and the versality of $\Omega$, after a suitable reparametrization $\nu \mapsto (\omega, \mu, \rho)$, the family $X$ takes the shape

$$X_{\omega,\mu,\rho} = [\omega + O(|y|,|z|)]\frac{\partial}{\partial x} + [\Omega(\mu)z + O(|y|,|z|^2)]\frac{\partial}{\partial z},$$

where $\mu = (\mu_1(\nu), \mu_2(\nu)) \in \mathbb{R}^2$ with $(\mu_1(\nu_0), \mu_2(\nu_0)) = (0, 0)$ and where $\Omega(\mu)$ is given by

$$\Omega(\mu) = \begin{pmatrix} 0 & -\lambda_0 - \mu_1 & 1 & 0 \\ \lambda_0 + \mu_1 & 0 & 0 & 1 \\ -\mu_2 & 0 & 0 & -\lambda_0 - \mu_1 \\ 0 & -\mu_2 & \lambda_0 + \mu_1 & 0 \end{pmatrix}. \quad (3.6)$$

---
[3]For a discussion on symplectic normal linearization, see [12, 26].

This matrix family is a *linear centralizer unfolding*[4] of $\Omega(0)$ in the matrix space $\mathfrak{sp}(4, \mathbb{R})$. This shows that the family $X = X_\nu$ always has two normal frequencies for all parameters $\nu$ sufficiently close to $\nu_0$, that is, $\omega^N(\nu) \in \mathbb{R}^2$. A geometric picture of the 'Cantor set' $\Gamma_{\tau,\gamma}(U)$ determined by the Diophantine conditions (3.2) is sketched in Figure 6, where we take $\nu = (\omega, \mu, \rho)$ and ignore the parameter $\rho$.[5] Variation in values of the parameter $\mu_2$ gives rise to a quasi-periodic Hamiltonian bifurcation. For this reason it is called the *detuning* (or distinguished) parameter of the bifurcation.

From the normal linear stability Theorem 3.3, we obtain the persistence of the invariant $m$-torus family that contains a normally $1 : -1$ resonant torus. As a corollary of Theorem 3.3, we have

**Theorem 3.5 (Persistence of Diophantine $m$-tori).** *Let $H = H_\nu$ be a $p$-parameter real-analytic family of integrable Hamiltonians given by (2.2). Suppose that*

- *The family $H$ is non-degenerate at the invariant torus $T_{\nu_0}$;*
- *The torus $T_{\nu_0}$ is normally generically $1 : -1$ resonant.*

*Then, for $\gamma$ sufficiently small and for any $p$-parameter real-analytic Hamiltonian family $\widetilde{H}$ on $(M, \sigma)$ sufficiently close to $H$ in the compact-open topology on complex analytic extensions, there exists a neighbourhood $\mathcal{U}$ of $\nu_0 \in \mathbb{R}^p$ and a map*

$$\Phi : M \times \mathcal{U} \to M \times \mathbb{R}^p,$$

*defined near the normally resonant tori $T_{\nu_0}$, such that,*

  i. *$\Phi$ is a $C^\infty$-smooth diffeomorphism onto its image and is a $C^\infty$-near the identity map;*
  ii. *The image $\Phi(V)$, where $V = \mathbb{T}^m \times \{(y,z) = (0,0)\} \times \Gamma_{\tau,\gamma}(\mathcal{U}')$, is a Cantor family of $\widetilde{H}$-invariant Diophantine tori, and the restriction of $\Phi$ to $V$ induces a conjugacy between $H$ and $\widetilde{H}$;*
  iii. *The restriction $\Phi|_V$ is symplectic and preserves the (symplectic) normal linear part $NX = \omega(\nu)\frac{\partial}{\partial x} + \Omega(\nu)z\frac{\partial}{\partial z}$ of the Hamiltonian vector field $X = X_H$ associated with $H$.*

---

[4]This is a linear versal unfolding with minimal number of parameters [3, 22].
[5]Note that internal as well as normal frequencies are independent of the parameter $\rho$.

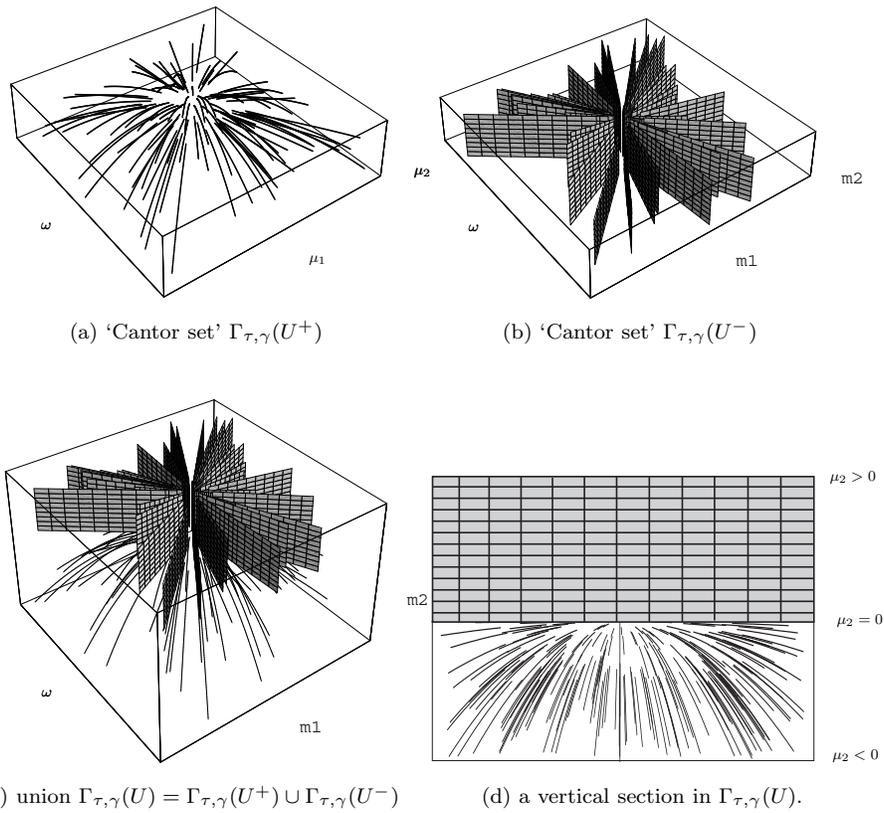

(a) 'Cantor set' $\Gamma_{\tau,\gamma}(U^+)$

(b) 'Cantor set' $\Gamma_{\tau,\gamma}(U^-)$

(c) union $\Gamma_{\tau,\gamma}(U) = \Gamma_{\tau,\gamma}(U^+) \cup \Gamma_{\tau,\gamma}(U^-)$

(d) a vertical section in $\Gamma_{\tau,\gamma}(U)$.

FIG 6. *Sketch of the 'Cantor sets' $\Gamma_{\tau,\gamma}(U^+)$ and $\Gamma_{\tau,\gamma}(U^-)$ corresponding to $\mu_2 \geq 0$ and $\mu_2 < 0$ respectively. The total 'Cantor set' $\Gamma_{\tau,\gamma}(U)$ is depicted in (c). The half planes in (b) and (c) give continua of invariant m-tori. In (d), a section of the 'Cantor set' $\Gamma_{\tau,\gamma}(U)$, along the $\mu_2$-axis, is singled out: the above grey region corresponds to a half plane given in (c)—a continuum of m-tori.*

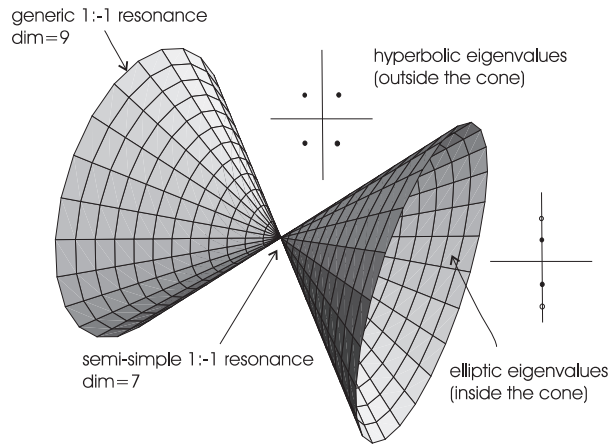

FIG 7. *Stratification of the 10-dimensional matrix space $\mathfrak{sp}(4,\mathbb{R})$. The cone represents the set of $1:-1$ resonant matrices. The subset of semi-simple matrices is given by the vertex of the cone, while the regular part of the cone contains the generic (or non-semisimple) ones.*

## 4. The quasi-periodic Hamiltonian Hopf bifurcation: persistence of the $(m+1)$- and $(m+2)$-tori

As mentioned before, see Section 2.2, the local torus foliation of the phase space $M$ near the normally $1 : -1$ resonant $m$-torus $T_{\nu_0}$ gives a stratification of a certain parameter space. This local stratification by tori is described by a piece of the swallowtail, see Figure 5(a). The smooth part of this surface corresponds to elliptic $(m+1)$-tori, while the open region above the surface (excluding the thread) corresponds to the Lagrangian $(m+2)$-tori. The persistence question of these tori can be answered by the 'standard' KAM theory [12, 26, 31, 32, 33, 35, 36].

### 4.1. Elliptic $(m+1)$-tori

To examine the persistence of the $(m + 1)$-tori, we first normalize the nearly-integrable perturbations of $H$ to obtain a proper integrable-approximation. This allows us to investigate the existence of invariant elliptic $(m + 1)$-tori of the perturbed Hamiltonian. To enable such a normalization, we need Theorem 3.5.

Let $\widetilde{H} = \widetilde{H}_\nu$ be a $p$-parameter real-analytic family of nearly-integrable Hamiltonian functions on $M = \mathbb{T}^m \times \mathbb{R}^m \times \mathbb{R}^4 = \{x, y, z\}$. We consider the quasi-periodic dynamics of the Hamiltonian $\widetilde{H}$, for $\widetilde{H}$ sufficiently close to the integrable family $H$ given by (2.2). First of all, by Theorem 3.5, the family $\widetilde{H}$ possesses a Cantor family of invariant $m$-tori. Secondly, the family $\widetilde{H}$ has the same normal linear behaviour as the integrable Hamiltonian $H$. Moreover, by the Inverse Function Theory and the versality of $\Omega$, a suitable reparametrization brings the family $\widetilde{H}$ into the form

$$\widetilde{H}_{\omega,\mu,\rho}(x,y,z) = \langle \omega, y \rangle + \frac{1}{2}\langle Jz, \Omega(\mu)z \rangle + \text{h.o.t.} , \qquad (4.1)$$

where the parameters $(\omega, \mu, \rho)$ are restricted to a 'Cantor set' and where the Floquet matrix $\Omega(\mu)$ is given by (3.6). To investigate the existence of the elliptic $(m+1)$-tori in the nearly-integrable family $\widetilde{H}$, we need to consider the higher order terms of the (normalized) Hamiltonian $\widetilde{H}$. The idea of normalization is to remove non-integrable terms from lower order terms by applying suitable coordinate transformations. These transformation are chosen to be the time-1 flows generated by certain Hamiltonian functions, for details see [7]. Such a normalization leads to the following:

**Theorem 4.1 (Quasi-periodic normal form [7]).** *Let $H = H_\nu$ be the real-analytic family of integrable Hamiltonians given by (2.2). Assume that assumptions of Theorem 3.3 are satisfied and that the parameter $\nu = (\omega, \mu, \rho) \in \Gamma_{\tau,\gamma}(U')$, see Section* 3.1, *where $U$ is a small neighbourhood of the fixed parameter $\nu_0$. Then, for any real-analytic Hamiltonian family $\widetilde{H} = \widetilde{H}_\nu$ sufficiently close to $H$ in the compact-open topology on complex analytic extensions, there exists a family of symplectic maps*

$$\Phi : \mathbb{T}^m \times \mathbb{R}^m \times \mathbb{R}^4 \times U \to \mathbb{T}^m \times \mathbb{R}^m \times \mathbb{R}^4 \times \mathbb{R}^p$$

*being real-analytic in $(x, y, z)$ and $C^\infty$-near-the-identity such that: the Hamiltonian $G = \widetilde{H} \circ \Phi$ is decomposed into the integrable part $G_{int}$ and a remainder $R$, where*

$$G_{int} = \langle \omega, y \rangle + (\lambda_0 + \mu_1)S + N + \mu_2 M + 2bM^2 + 2c_1 SM + c_2 S^2 ,$$

with $S = z_1 z_4 - z_2 z_3$, $N = \frac{1}{2}(z_3^2 + z_4^2)$ and $M = \frac{1}{2}(z_1^2 + z_2^2)$. The remainder $R$ satisfies

$$\frac{\partial^{q+|p|+|k|} R}{\partial^q \mu_2 \partial^p y \partial^k z}(x, 0, 0, \omega, \mu_1, 0) = 0,$$

for all indices $(q, p, k) \in \mathbb{N}_0^3$ with $2q + 4|p| + |k| \leq 4$.

For a similar non-linear normal form theory for the non-conservative setting see [5, 12] and for the Hamiltonian setting see [23]. A special case of Theorem 4.1 with $m = 1$ is extensively considered in [34].

Now the integrable truncation $G_{int}$ contains a family of elliptic $(m+1)$-tori determined by the cubic equation

$$S^2 - 4bM^3 - 4\mu_2 M^2 - 4c_1 SM^2 = 0, \qquad (4.2)$$

where $M > 0$, compare with [17, 24, 29]. By a rescaling argument, the remainder $R$ can be considered as perturbation of the integrable $G_{int}$. At this point, our present problem is cast into the form treated in Theorem 2.6 of [11] from which we conclude that most elliptic $(m+1)$-tori of the integrable Hamiltonian $G_{int}$ survive the perturbations by $R$, and are only slightly deformed.

**Theorem 4.2 (Persistence of Diophantine elliptic $(m+1)$-tori).** *Let $H = H_\nu$ be the real-analytic family of integrable Hamiltonian functions given by (2.2) such that assumptions from Theorem 3.5 are satisfied. Then, for any real-analytic Hamiltonian family $\widetilde{H} = \widetilde{H}_\nu$ sufficiently close to the family $H$ in the compact-open topology, there exists a map*

$$\Phi : \mathbb{T}^m \times \mathbb{R}^m \times \mathbb{R}^4 \times U \to \mathbb{T}^{m+1} \times \mathbb{R}^{m+1} \times \mathbb{R}^2 \times \mathbb{R}^p$$

*defined near the normally resonant torus $T_{\nu_0}$ such that: $\Phi$ is a $C^\infty$-near-the-identity diffeomorphism onto its image and the Hamiltonian $\widetilde{H} \circ \Phi^{-1}$ has a Cantor family of invariant elliptic $(m+1)$-tori.*

### 4.2. Lagrangian $(m+2)$-tori

In this section, we investigate the persistence of the Lagrangian $(m + 2)$-tori of the integrable Hamiltonian $H$ of the form (2.2). These tori are located in the open region (excluding the thread) above the swallowtail surface, see Figure 5(a). To this end, we apply the classical KAM theory [2, 27, 33] and its global version [6]. The Kolmogorov non-degeneracy condition — required for the KAM theory — on the Lagrangian tori near the resonant torus is guaranteed by the non-triviality of the monodromy, compare with [19, 38, 45].

Reconsider the integrable (i.e., $\mathbb{T}^m$-symmetric) Hamiltonian function $H$ of the form

$$H(x, y, z, \nu) = \langle \omega(\nu), y \rangle + \frac{1}{2} \langle Jz, \Omega(\nu)z \rangle + F(y, z, \nu), \qquad (4.3)$$

on the space $M = \mathbb{T}^m \times \mathbb{R}^m \times \mathbb{R}^4 = \{x, y, z\}$, where $F$ denotes the higher order terms. The invariant torus $T_{\nu_0}$ given by $(y, z, \nu) = (0, 0, \nu_0)$ is generically $1 : -1$ resonant, see Section 2.2. We require that the Hessian of the higher order term $F$ with respect to the variable $y$ is non-vanishing at the torus $T_{\nu_0}$.

As before we may assume that $H$ is invariant under the free $\mathbb{S}^1$-action generated by the semisimple part of the polynomial $\frac{1}{2}\langle Jz, \Omega(\nu)z \rangle$, see Section 2.2. By this invariance, the Hamiltonian $H$ is Liouville-integrable with the $(m+2)$ first integrals $\mathcal{EM} = (H, S, y)$, where $S = z_1 z_4 - z_2 z_3$. The Lagrangian $(m+2)$-tori of the integrable Hamiltonian $H$ are the regular fibers of the energy-momentum map $\mathcal{EM}$. Our concern is with the persistence of these Lagrangian tori near the thread, see Figure 5(a), when $H$ is perturbed into a nearly-integrable Hamiltonian family $\widetilde{H}$. In view of the classical KAM theory [2, 27], most of these Lagrangian tori survive the perturbation in Whitney-smooth Cantor families. Here we have to require the Kolmogorov non-degeneracy, which near the thread is a consequence of the non-trivial monodromy. Indeed, non-degenacy follows by an application of [19, 38, 40, 45] to the reduced two-degrees-of-freedom system $H_\alpha$ and the assumption that $\det \frac{\partial^2 F}{\partial y^2} \neq 0$ at the resonant torus $T_{\nu_0}$. We first conclude that, for sufficiently small perturbation, the Lagrangian tori survive in a Whitney-smooth Cantor family of positive measure [2, 33, 35]. Secondly, the corresponding KAM-conjugacies, which are only defined on locally trivial sub-bundles, can be glued together to provide a globally Whitney-smooth conjugacy from the integrable to the nearly-integrable Cantor torus family [6].

**Theorem 4.3 (Persistence of Diophantine Lagrangian $(m+2)$-tori).** *Let $H = H_\nu$ be the real-analytic family of Hamiltonians given by (4.3). Suppose that the Hessian of the higher order term $F$, see (4.3), with respect to $y$ is non-zero at the resonant torus $T_{\nu_0}$. Then, there exists a neighbourhood $U \subset \mathbb{R}^p$ of $\nu_0$ such that for any real-analytic Hamiltonian $\widetilde{H}$ sufficiently close to $H$ in the compact-open topology on complex analytic extensions, the following holds: the perturbed Hamiltonian $\widetilde{H}$ has a Cantor family of invariant $(m+2)$-tori; this family is a $C^\infty$-near-the-identity diffeomorphic image of $T = \bigcup_\nu T_\nu$, where $\nu$ is restricted to a 'Cantor set' (determined by the Diophantine conditions on the internal frequencies); in these tori, the diffeomorphism conjugates $H$ and $\widetilde{H}$.*

## 5. Concluding remarks

We considered a family of $\mathbb{T}^m$-symmetric Hamiltonians on $M = \mathbb{T}^m \times \mathbb{R}^m \times \mathbb{R}^4$ that has a normally $1:-1$ resonant torus. As the parameter varies, the torus changes from normally hyperbolic into normally elliptic. This generically gives rise to a quasi-periodic Hamiltonian Hopf bifurcation. Near the normally resonant torus the phase space $M$ is foliated by hyperbolic and elliptic $m$-tori, by elliptic $(m+1)$-tori and by Lagrangian $(m+2)$-tori. This singular torus foliation gives a stratification in a suitable parameter space: the strata are determined by the dimension of the tori. The local geometry of this stratification is a piece of swallowtail catastrophe set: the $m$-tori are located at the 1-dimensional part of the surface, $(m+1)$-tori at the regular part of the surface and $(m+2)$-tori at the open region above the surface, see Figure 5(a). By KAM-theory [2, 10, 12, 26, 27, 33, 35, 36], these tori survive in Whitney-smooth Cantor families, under small nearly-integrable perturbations, compare with Figure 5(b). We remark that this quasi-periodic stability still holds for the case where the frequencies of the oscillator are kept constant, see [7] for details. In view of the global KAM theory [6], the non-trivial monodromy in the integrable Lagrangian torus bundle can be extended to the nearly-integrable case. This may be of interest for semi-classical quantum mechanics, compare with [14, 15, 18, 20, 21, 41]. An example of the quasi-periodic Hamiltonian Hopf bifurcation is the Lagrange top (near gyroscopic stabilization) coupled to a quasi-periodic oscillator, compare with Section 1.

Concerning the parameter domains regarding the Diophantine tori of the different dimensions $m, m+1$ and $m+2$, we expect that these are attached to one another in a Whitney-smooth way, as suggested by the integrable approximation. Following [8] we speak of a Cantor stratification. Here the stratum of the $(m+1)$-tori consists of density points of $(m+2)$-quasiperiodicity of $(2m+4)$-dimensional Hausdorff measure. Similarly, the stratum of the $m$-tori consists of density points of $(m+1)$-quasiperiodicity of $(2m+2)$-dimensional Hausdorff measure. Also compare with [11].

In a general case where the perturbation destroys the whole $\mathbb{T}^m$-symmetry, the Cantor families of stable KAM-tori always contains continua of tori (the projection of these continua into the parameter space has no Cantor gaps). We refer to these projections as continuous structures in the Cantor family of surviving tori. For special perturbations where a partial symmetry still remains (e.g., perturbations that are independent of certain angle variables), we expect extra continuous structures, for details see [7].

In this paper, we addressed the supercritical case of the quasi-periodic Hamiltonian Hopf bifurcation. We expect similar persistence results for the subcritical case. This is important for a better understanding of the hydrogen atom in crossed electric and magnetic fields [20, 21], where the subcritical bifurcation occurs. The singular torus foliation for this case is described by the tail of the swallowtail surface.[6] We expect that our approach for the supercritical case also works for the subcritical case.

## Acknowledgments

The authors are grateful to Floris Takens for his valuable comments.

---

[6]Complementary to the part of Figure 5(a).